\documentclass[11pt,a4paper]{article}
\usepackage[margin=28mm]{geometry}
\usepackage[T1]{fontenc}
\usepackage{lmodern,microtype}
\usepackage{amsmath,amssymb,amsthm,mathtools}
\usepackage{enumitem,needspace,color,graphicx}
\usepackage[hidelinks]{hyperref}
\hypersetup{pdftitle={Optimal no-feedback card guessing after a fixed number of riffle shuffles},
 pdfauthor={Tipaluck Krityakierne and Thotsaporn Aek Thanatipanonda}}
\numberwithin{equation}{section}
\newtheorem{theorem}{Theorem}[section]
\newtheorem{lemma}[theorem]{Lemma}
\newtheorem{proposition}[theorem]{Proposition}
\newtheorem{corollary}[theorem]{Corollary}
\theoremstyle{remark}
\newtheorem{remark}[theorem]{Remark}
\newcommand{\Prob}{\mathbb P}
\newcommand{\E}{\mathbb E}
\newcommand{\Bin}{\operatorname{Bin}}
\newcommand{\coef}[1]{[x^{#1}]}
\newcommand{\Z}{\mathbb Z}
\newcommand{\ind}{\mathbf 1}
\setlist[enumerate]{label=\textup{(\roman*)},leftmargin=*,itemsep=3pt}
\title{Optimal no-feedback card guessing after a fixed number of riffle shuffles}
\author{Tipaluck Krityakierne\thanks{Department of Mathematics, Faculty of Science,
Mahidol University, Bangkok 10400, Thailand; and Centre of Excellence in Mathematics,
MHESI, Bangkok 10400, Thailand. Email: \texttt{tipaluck.kri@mahidol.ac.th}.}
\and Thotsaporn Aek Thanatipanonda\thanks{Science Division, Mahidol University
International College, Nakhon Pathom 73170, Thailand.
Email: \texttt{thotsaporn@gmail.com}.}}
\date{}

\begin{document}
\maketitle
\begin{abstract}
An ordered deck of $n$ cards is subjected to $k$ independent
Gilbert--Shannon--Reeds riffle shuffles. A player then guesses each card
without feedback, with repeated guesses allowed. For every fixed $k\ge1$,
we prove that, once $n$ is sufficiently large, an optimal guess at position
$i$ in the top half of the deck is $\lfloor i/2^k\rfloor+1$; the bottom-half
strategy is obtained by reflection. We also show that the optimal expected
number of correct guesses is $2\sqrt{n/((2^k-1)\pi)}+O_k(1)$.
Our starting point is the result of our earlier work that the label at
each position is distributed as an equally weighted mixture of $2^k$
components, each a shifted sum of two independent binomial random
variables (a Poisson-binomial distribution). The proof of the optimal strategy combines three
ingredients: a uniform concentration estimate that separates the mixture
components, exact binomial inequalities that control the gaps between
competing peaks, and a coefficient comparison that resolves modal ties.
The expected-score asymptotic is proved independently of the strategy,
using only the peak heights of the mixture components.
\end{abstract}


\section{Introduction and main results}\label{sec:introduction}

Consider a deck of $n$ distinct cards, initially ordered $1,2,\ldots,n$
from top to bottom. The deck is shuffled $k$ times according to the
Gilbert--Shannon--Reeds model. In one shuffle, the deck is cut into a top packet of $c$ cards and a
bottom packet of $n-c$ cards, where $c$ has the binomial distribution
$\Bin(n,1/2)$; the two packets are then interleaved while preserving
their internal orders: at each step, the next card is taken from a packet
chosen with probability proportional to the number of cards remaining in it. Successive shuffles are independent.
This is the classical riffle-shuffle model; see Bayer and
Diaconis~\cite{BD92}.

In the no-feedback guessing game, the player chooses one card label for
each position without learning the identities of the cards or the
correctness of any previous guesses. The player may repeat a label at
different positions. Thus a strategy is simply a sequence
$G=(g_1,\ldots,g_n)\in\{1,\ldots,n\}^n$.
Following~\cite{KT25}, let $J_i$ be the label at position $i$ after the
$k$ shuffles and let
\[
 a^{(k)}_{i,j}=\Prob(J_i=j),\qquad
 P^{(k)}=(a^{(k)}_{i,j})_{1\le i,j\le n}.
\]
The number of correct guesses under $G$ is
$X^G=\sum_{i=1}^n\ind_{\{J_i=g_i\}}$. Consequently,
\begin{equation}\label{eq:objective}
 \E X^G=\sum_{i=1}^n a^{(k)}_{i,g_i},\qquad
 \E X^{G^*}=\sum_{i=1}^n\max_{1\le j\le n}a^{(k)}_{i,j},
\end{equation}
where $G^*$ denotes an optimal strategy. The optimization therefore
depends only on the marginal distribution at each position.

For a single dovetail (that is, riffle) shuffle, Ciucu~\cite{Ciucu98}
showed that an optimal guessing sequence begins $1,2,2,3,3,\ldots$ in the
top half of the deck and is reflected in the bottom half, and that the
optimal expected score is $2\sqrt{n/\pi}+O(1)$. He also analyzed optimal
guessing when the number of shuffles grows logarithmically with the deck
size. Subsequent work investigated the one-shuffle
score in greater detail. Krityakierne et al.~\cite{KSTY23} derived moment
formulas, while Kuba and Panholzer~\cite{KP24} established a limit law
and convergence of the integer moments. The complete-feedback version of
the game after one riffle shuffle was treated with generating functions
in~\cite{KT23}.

Returning to the no-feedback game, our earlier paper~\cite{KT25}
provides an explicit generating-function representation of every marginal $a^{(k)}_{i,j}$ for arbitrary $n$ and $k$.
For each fixed position $i$, the marginal distribution of $J_i$ is an
equally weighted mixture of $2^k$ distributions, each given by a shifted sum of two independent binomial random variables.
That representation suggests the optimal strategy when $k$ is fixed and
$n$ tends to infinity, but it does not by itself identify the label at
which the mixture is largest. 

The present paper establishes the exact optimal strategy for all
sufficiently large $n$, together with the asymptotic optimal expected
score. The threshold $n(k)$ in the first theorem
is uniform over all positions in the deck.

\begin{theorem}[Optimal strategy for fixed $k$]\label{thm:strategy}
For each fixed integer $k\ge1$, there is an integer $n(k)$ such that,
for every $n\ge n(k)$, the strategy
\begin{equation}\label{eq:strategy}
 j_i^*=\begin{cases}
 \displaystyle\left\lfloor\frac{i}{2^k}\right\rfloor+1,
 &\displaystyle 1\le i\le\left\lfloor\frac{n+1}{2}\right\rfloor,\\[7pt]
 \displaystyle n-\left\lfloor\frac{n+1-i}{2^k}\right\rfloor,
 &\displaystyle\left\lfloor\frac{n+1}{2}\right\rfloor<i\le n
 \end{cases} 
\end{equation} 
is optimal. Equivalently, $a^{(k)}_{i,j_i^*}=\max_j a^{(k)}_{i,j}$ for every
$i$. If $n$ is odd, both $j_{(n+1)/2}^*$ and its reflection
$n+1-j_{(n+1)/2}^*$ are optimal at the middle position.
\end{theorem}

Thus the top-half sequence consists of $2^k-1$ initial guesses of $1$,
followed by blocks of $2^k$ guesses of $2,3,\ldots$, truncated at the
middle of the deck. The bottom-half sequence is its reflection.
The theorem identifies one optimal sequence; additional
optimal choices can occur. For $k\ge2$, we determine all such choices
for sufficiently large $n$ in Corollary~\ref{cor:uniqueness}.

\begin{theorem}[Optimal expected score]\label{thm:expectation}
For every fixed integer $k\ge1$,
\begin{equation}\label{eq:expectation}
 \E X^{G^*}=\frac{2\sqrt n}{\sqrt{(2^k-1)\pi}}+O_k(1)
 \qquad\text{as }n\longrightarrow\infty.
\end{equation}
\end{theorem}

Theorem~\ref{thm:strategy} settles the fixed-$k$ optimal-strategy question
posed in Conjecture 4.1 of~\cite{KT25}, with a small adjustment
to the bottom-half formula at divisibility ties. Theorem~\ref{thm:expectation} extends the
one-shuffle leading term to every fixed number of shuffles. Its proof
uses a uniform estimate for the row maxima and does not depend on
Theorem~\ref{thm:strategy} or the resolution of ties.

The argument rests on three observations, the first two of which are illustrated in Figure~\ref{fig:mixture}. The component means are spaced
by $(n-1)/2^k$, so at any given label all but the closest component have
exponentially small probability. In the top half, the first component
has the largest peak height for sufficiently large $n$. Finally, every
endpoint gap, that is, the gap between the peak height of
the first or last component and its value at a nonmodal label, or between
these two peak heights (Section~\ref{sec:representation}),
has a uniform polynomial
lower bound obtained directly from binomial ratios. These facts reduce
the exact optimization to a small number of tie comparisons.

Section~\ref{sec:representation} recalls the representation from~\cite{KT25}.
Section~\ref{sec:estimates} establishes the required estimates.
Sections~\ref{sec:strategy} and~\ref{sec:expectation} prove the two main
theorems. Section~\ref{sec:discussion} records the scope of the results
and further questions. Throughout, $k$ is fixed, and constants in
$O_k(\cdot)$ may depend on $k$.

\section{The marginal generating functions}\label{sec:representation}

We retain the notation and component indexing of~\cite{KT25}. Set
$m=2^k$ and $p=1/m$; since $m$ is determined by $k$, we write $O_m(\cdot)$
and $O_k(\cdot)$ interchangeably. For fixed $n$ and $i$, Theorem 3.3 of that paper gives the probability
generating function $F_{L_t}$ of the $t$th mixture component:
\begin{equation}\label{eq:component}
 F_{L_t}(x)
 =x\left(1-\frac{t}{m}+\frac{t}{m}x\right)^{n-i}
       \left(1-\frac{t+1}{m}+\frac{t+1}{m}x\right)^{i-1},
 \quad 0\le t<m.
\end{equation}
Here $L_t$ is the binary-list index defined by the bisection construction
in~\cite{KT25}. We use its order $t=0,\ldots,m-1$ without revisiting
that construction. With $p_{L_t}(j)=\coef{j}F_{L_t}(x)$, 
the probability mass function of
the $t$th component, 
Corollary 3.5
of~\cite{KT25} states that
\begin{equation}\label{eq:mixture}
 G_i(x):=\E x^{J_i}=\frac1m\sum_{t=0}^{m-1}F_{L_t}(x),
 \qquad
 a^{(k)}_{i,j}=\frac1m\sum_{t=0}^{m-1}p_{L_t}(j).
\end{equation}
The dependence of $F_{L_t}$ and $p_{L_t}$ on $n,i,k$ is suppressed,
as in~\cite{KT25}. Likewise, the means $\mu_t$ in~\eqref{eq:means} depend
on the row, and $\beta_N,H_N$ in~\eqref{eq:binomialnotation} depend
on~$p$. Unless stated otherwise, $\max_j$ denotes the maximum over the
labels $1\le j\le n$.

The factorization in~\eqref{eq:component} identifies $p_{L_t}$ as the
distribution of
\begin{equation}\label{eq:bernoulli}
 1+\Bin\left(n-i,\frac{t}{m}\right)
   +\Bin\left(i-1,\frac{t+1}{m}\right),
\end{equation}
with independent binomial summands. For the analysis of a fixed row $i$,
write
\begin{equation}\label{eq:auxiliary}
 N=i-1,\qquad M=n-i,\qquad d=n-1=N+M,\qquad v=p(1-p).
\end{equation}
The mean of the $t$th component is
\begin{equation}\label{eq:means}
 \mu_t=1+\frac{dt+N}{m}\quad(0\le t\le m-1),\qquad
 \mu_{t+1}-\mu_t=\frac d m\quad(0\le t\le m-2).
\end{equation}

\begin{figure}[t]
\centering
\includegraphics[width=0.9\textwidth]{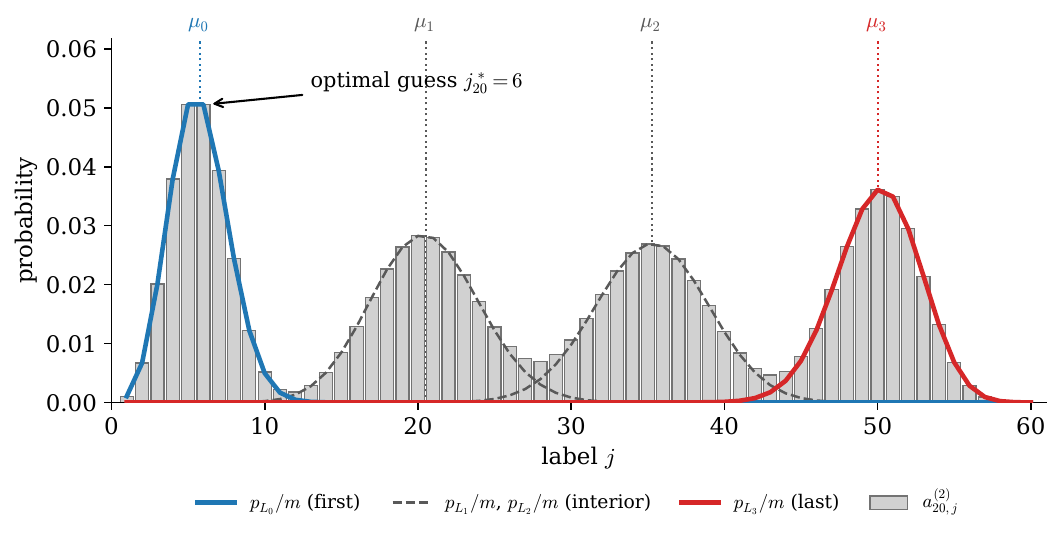}
\caption{The row $a^{(2)}_{20,j}$ for $n=60$ and $k=2$ (bars), decomposed
by~\eqref{eq:mixture} into its $m=4$ components $p_{L_t}/m$ (curves),
with the means $\mu_t$ of~\eqref{eq:means} marked. The first component
($t=0$) has the tallest peak and contains the optimal guess
$j^*_{20}=\lfloor 20/4\rfloor+1=6$. Since $20/4$ is an integer, the first
component has two modes, $5$ and $6$; the other components break this
tie in favour of $6$ (Case~2.1 in Section~\ref{sec:strategy}), although only
by $a^{(2)}_{20,6}-a^{(2)}_{20,5}\approx 9.3\times10^{-7}$, far too little to be
visible in the plot.}
\label{fig:mixture}
\end{figure}

Figure~\ref{fig:mixture} shows a typical row. The components are
separated by $d/m$, so at each label essentially only the component with
the nearest mean contributes (Lemma~\ref{lem:tails}). For sufficiently
large $n$, in the top half the first component $t=0$ has a maximal peak
height, and its mode $\lfloor i/m\rfloor+1$ (the larger one, when there
are two) is an optimal guess. The only possible competitor is then the
last component $t=m-1$, whose peak height is compared with that of the first component
in Lemma~\ref{lem:binomial}(ii); the interior components are strictly
lower (Lemma~\ref{lem:interior}).

We use two auxiliary binomial quantities, defined for every integer
$N\ge0$ (not only for the value $N=i-1$ of~\eqref{eq:auxiliary}):
\begin{equation}\label{eq:binomialnotation}
 \beta_N(r)=\Prob(\Bin(N,p)=r),\qquad
 H_N=\max_{r\in\Z}\beta_N(r),
\end{equation}
where $\beta_N(r)=0$ outside $0\le r\le N$. In particular, $H_0=1$.
With $N$ and $M$ as in~\eqref{eq:auxiliary}, the first and
last components $t=0$ and $t=m-1$, which we call the \emph{endpoint
components},
satisfy
\begin{equation}\label{eq:endpoints}
 p_{L_0}(j)=\beta_N(j-1),\qquad
 p_{L_{m-1}}(j)=\beta_M(n-j).
\end{equation}
Their peak heights are $H_N$ and $H_M$, respectively.
We call $H_N$ and $H_M$ the \emph{endpoint peaks}. An
\emph{endpoint gap} is a positive difference of one of two types:
$H_N-\beta_N(r)$ or $H_M-\beta_M(r)$ with $r$ nonmodal (a peak minus a
nonmodal value of the same endpoint component), or $H_N-H_M$ when
$H_N>H_M$ (the difference of the two endpoint peaks).
Lemma~\ref{lem:binomial} bounds them below by~\eqref{eq:pointgap}
and~\eqref{eq:heightgap}, respectively.

The reflection symmetry of the marginal matrix is
\begin{equation}\label{eq:symmetry}
 a^{(k)}_{i,j}=a^{(k)}_{n+1-i,n+1-j}.
\end{equation}
It follows either from the one-shuffle symmetry in~\cite{KT25} and
matrix powers, or directly from~\eqref{eq:component} by replacing
$i$ with $n+1-i$ and $t$ with $m-1-t$.
Hence it suffices to analyze $i\le\lfloor(n+1)/2\rfloor$, where $N\le M$.
We call $t=1,\ldots,m-2$ the interior components. There are no interior
components when $k=1$.

\section{Uniform estimates for the components}\label{sec:estimates}

\subsection{Exact binomial gaps}

The following elementary estimates provide lower bounds for the gaps
between binomial peak heights, uniform over all positions, including
those close to the ends of the deck.
The polynomial order is sufficient for comparison with exponential tails.

\begin{lemma}\label{lem:binomial}
Let $m\ge2$ be any integer (we apply the lemma with $m=2^k$), let
$p=1/m$ and $v=p(1-p)$, and define $\beta_N,H_N$
by~\eqref{eq:binomialnotation}. In this lemma $N$ and $M$ denote arbitrary
nonnegative integers. Then:
\begin{enumerate}
\item The modes of $\beta_N$ are $\lfloor(N+1)p\rfloor$ and, when
$(N+1)p$ is an integer, also $(N+1)p-1$. For every nonmodal $r\in\Z$,
\begin{equation}\label{eq:pointgap}
 H_N-\beta_N(r)\ge\frac1{m(N+1)^2}.
\end{equation}
\item $H_{N+1}\le H_N$, with equality if and only if $m\mid (N+1)$.
For integers $0\le N<M$,
\begin{equation}\label{eq:equality}
 H_N=H_M\quad\Longleftrightarrow\quad
 M=N+1\text{ and }m\mid (N+1).
\end{equation}
If $H_N>H_M$, then
\begin{equation}\label{eq:heightgap}
 H_N-H_M\ge\frac1{m^2(M+1)^2}.
\end{equation}
\item As $N\to\infty$,
\begin{equation}\label{eq:stirling}
 H_N=\frac1{\sqrt{2\pi vN}}+O_m(N^{-3/2}).
\end{equation}
\end{enumerate}
\end{lemma}

\begin{proof}
The ratio
\[
 \frac{\beta_N(r+1)}{\beta_N(r)}
 =\frac{(N-r)p}{(r+1)(1-p)},\qquad 0\le r<N,
\]
gives the modes and unimodality. Let $r$ be a mode of $\beta_N$ with a
nonmodal neighbor $r+1$ or $r-1$. The differences $H_N-\beta_N(r+1)$ and
$H_N-\beta_N(r-1)$ are, respectively,
\[
 H_N\frac{r+1-(N+1)p}{(r+1)(1-p)},\qquad
 H_N\frac{(N+1)p-r}{(N-r+1)p}.
\]
The numerator corresponding to the nonmodal neighbor is positive, and
both numerators lie in $\frac1m\Z$ because $(N+1)p=(N+1)/m$; hence that
numerator is at least $1/m$. Each denominator is at most $N+1$, and
$H_N\ge1/(N+1)$ because the $N+1$ masses of $\beta_N$ sum to $1$.
This proves~\eqref{eq:pointgap}
at the adjacent nonmodal values; unimodality gives the result everywhere
else within the support. Outside the support it follows directly from
$H_N\ge1/(N+1)$. This also covers $N=0$.

The binomial recursion gives
\begin{equation}\label{eq:binrecursion}
 \beta_{N+1}(r)=(1-p)\beta_N(r)+p\beta_N(r-1).
\end{equation}
Thus $H_{N+1}\le H_N$. Equality holds precisely when both $r$ and $r-1$
are modes of $\beta_N$ for some $r$, that is, when $m\mid (N+1)$.
If the inequality is strict, take $r$ maximizing the left-hand side
of~\eqref{eq:binrecursion}. At least one of $r,r-1$ is nonmodal for
$\beta_N$, and its weight in~\eqref{eq:binrecursion} is at least $p$
(since $p\le1/2$). By~\eqref{eq:pointgap}, therefore,
\[
 H_N-H_{N+1}\ge\frac1{m^2(N+1)^2}.
\]
Two consecutive steps cannot both give equality, since $m\ge2$; with
monotonicity this proves~\eqref{eq:equality}. For~\eqref{eq:heightgap},
if $H_N>H_M$ then some step $N\le N'<M$ is strict, and monotonicity gives
\[
 H_N-H_M\ge H_{N'}-H_{N'+1}\ge\frac1{m^2(N'+1)^2}
 \ge\frac1{m^2(M+1)^2}.
\]

Finally, a binomial mode differs from $Np$ by less than $1$.
Stirling's formula (equivalently, the de Moivre--Laplace local limit
theorem), applied at such an integer $r$, yields
\[
 \beta_N(r)=\frac1{\sqrt{2\pi Nv}}
              \bigl(1+O_m(N^{-1})\bigr).
\]
The bound is uniform over the two possible modes, and proves~\eqref{eq:stirling}.
\end{proof}

\subsection{Interior peaks and separated tails}

\begin{lemma}[Interior peak bound]\label{lem:interior}
For fixed $m=2^k\ge4$, there are constants $\kappa_m>0$ and $n_0(m)$ such that,
for every $n\ge n_0(m)$, every $1\le i\le\lfloor(n+1)/2\rfloor$, and every
$1\le t\le m-2$,
\begin{equation}\label{eq:interiorgap}
 \max_j p_{L_t}(j)\le H_{i-1}-\kappa_m n^{-1/2}.
\end{equation}

\end{lemma}

\begin{figure}[t]
\centering
\includegraphics[width=0.9\textwidth]{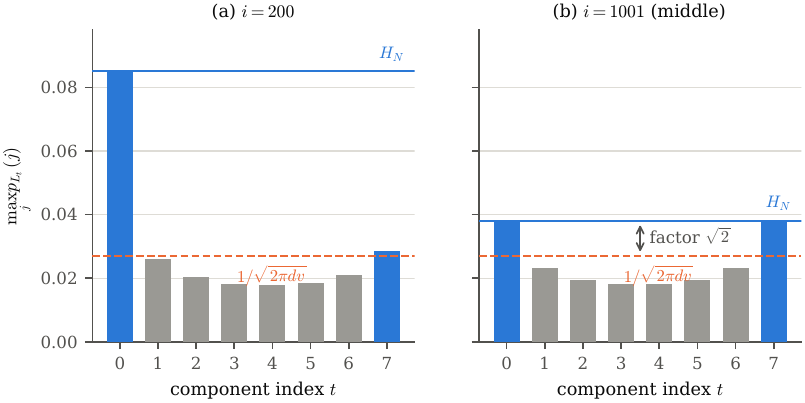}
\caption{Peak heights $\max_j p_{L_t}(j)$ of the $m=8$
components for $k=3$, $n=2001$: the endpoint peaks $H_N$, $H_M$ (blue)
and the interior peaks (grey). Every interior peak lies below the leading
term $1/\sqrt{2\pi dv}$ of~\eqref{eq:fourier} (dashed). (a) At a typical
top-half position, $H_N$ stands far above all other peaks. (b) At the
middle position, $H_N=H_M$ and the margin is smallest:
$H_N\approx(\pi vd)^{-1/2}$, a factor $\sqrt2$ above the dashed line.}
\label{fig:peaks}
\end{figure}

By~\eqref{eq:endpoints}, $H_{i-1}$ is the peak height of the first
component. Thus~\eqref{eq:interiorgap} says that every interior peak lies
below the first-component peak by at least $\kappa_mn^{-1/2}$, uniformly over
all positions in the top half. Figure~\ref{fig:peaks} shows
this for $k=3$: the gap is large at a typical position and smallest at the
middle of the deck.

\begin{proof}
In an interior component, all $d=n-1$ Bernoulli parameters in
\eqref{eq:bernoulli} belong to $[p,1-p]$, so $q(1-q)\ge v$ for each such
parameter $q$. The characteristic function
$\phi_q(\theta)=1-q+qe^{\mathrm i\theta}$ of a Bernoulli$(q)$ variable
satisfies
\[
 |\phi_q(\theta)|^2
 =1-4q(1-q)\sin^2(\theta/2)
 \le\exp\{-4q(1-q)\sin^2(\theta/2)\}
 \le\exp\{-4v\sin^2(\theta/2)\}.
\]
The characteristic function of the component is
\[
 \phi(\theta)=\sum_j p_{L_t}(j)e^{\mathrm ij\theta}
 =e^{\mathrm i\theta}\,\phi_{t/m}(\theta)^{n-i}\,\phi_{(t+1)/m}(\theta)^{i-1},
\]
that is, $e^{\mathrm i\theta}$ times a product of $d$ factors of the form
$\phi_q(\theta)$. Taking square
roots and multiplying gives
$|\phi(\theta)|\le\exp\{-2dv\sin^2(\theta/2)\}$. By Fourier inversion,
$p_{L_t}(j)=\frac1{2\pi}\int_{-\pi}^{\pi}\phi(\theta)e^{-\mathrm ij\theta}\,d\theta$
for every $j\in\Z$. Therefore
\begin{align}
 \max_j p_{L_t}(j)
 &\le\frac1{2\pi}\int_{-\pi}^{\pi}
       e^{-2dv\sin^2(\theta/2)}\,d\theta\notag\\
 &=\frac1{\sqrt{2\pi dv}}+O_m(d^{-3/2}).\label{eq:fourier}
\end{align}

To justify the last equality, compare the integrand with the Gaussian
$e^{-dv\theta^2/2}$, whose integral over $\mathbb R$, divided by $2\pi$,
is exactly $(2\pi dv)^{-1/2}$. Two error terms remain.
\begin{itemize}[leftmargin=*]
\item \emph{Gaussian tail.} The part of the Gaussian integral outside
$[-\pi,\pi]$ is $O(e^{-\pi^2dv/2})$, which is exponentially small
in $d$.
\item \emph{Difference of integrands on $[-\pi,\pi]$.} From
$x-x^3/6\le\sin x\le x$ for $0\le x\le\pi/2$ we get
$\theta^2/4-\theta^4/48\le\sin^2(\theta/2)\le\theta^2/4$, so the two
exponents $2dv\sin^2(\theta/2)$ and $dv\theta^2/2$ differ by at most
$dv\theta^4/24$. Also $\sin^2(\theta/2)\ge\theta^2/\pi^2$ on
$[-\pi,\pi]$. Using $|e^{-a}-e^{-b}|\le|a-b|\,e^{-\min(a,b)}$,
\[
 \bigl|e^{-2dv\sin^2(\theta/2)}-e^{-dv\theta^2/2}\bigr|
 \le\frac{dv\,\theta^4}{24}\,e^{-2dv\theta^2/\pi^2},
\]
and integrating over $\mathbb R$ gives
$\frac{dv}{24}\cdot\frac{3\sqrt\pi}{4}\bigl(\frac{2dv}{\pi^2}\bigr)^{-5/2}
=O_m(d^{-3/2})$.
\end{itemize}
Together these prove~\eqref{eq:fourier}.

To deduce~\eqref{eq:interiorgap}, we compare~\eqref{eq:fourier}
with the first-component peak $H_N$. Since $N=i-1\le\lfloor d/2\rfloor$, the monotonicity of
Lemma~\ref{lem:binomial}(ii) and the asymptotic formula
of Lemma~\ref{lem:binomial}(iii) give
\begin{equation}\label{eq:firstlower}
 H_N\ge H_{\lfloor d/2\rfloor}
 =\frac1{\sqrt{\pi vd}}+O_m(d^{-3/2}),
\end{equation}
where replacing $\lfloor d/2\rfloor$ by $d/2$ contributes a further
$O_m(d^{-3/2})$.
The leading constants in~\eqref{eq:fourier} and~\eqref{eq:firstlower}
differ by a factor of $\sqrt2$ (Figure~\ref{fig:peaks}(b)), so
$H_N-\max_j p_{L_t}(j)\ge(1-2^{-1/2})(\pi vd)^{-1/2}+O_m(d^{-3/2})$,
which is at least $\kappa_m n^{-1/2}$ for all sufficiently large $n$. Both estimates
are independent of $i$, proving the required uniformity.

\end{proof}

\begin{lemma}[Uniform separation estimate]\label{lem:tails}
For $n\ge2$, define
\begin{equation}\label{eq:epsilon}
 \varepsilon_n=2(m-1)\exp\left(-\frac{n-1}{2m^2}\right).
\end{equation}
For any position $i$ and label $j$, choose an index $s$ whose 
mean $\mu_s$ is closest to $j$. 
Then
\begin{equation}\label{eq:tails}
 0\le\sum_{t=0}^{m-1}p_{L_t}(j)-p_{L_s}(j)\le\varepsilon_n.
\end{equation}
\end{lemma}

In words, at every label $j$ the row value $a^{(k)}_{i,j}$ agrees with
$p_{L_s}(j)/m$, the weighted component whose mean is nearest, up to the
exponentially small error $\varepsilon_n/m$ (see Figure~\ref{fig:mixture}).

\begin{proof}
By~\eqref{eq:means}, the means are separated by $d/m$, where $d=n-1$.
If $t\ne s$, the choice of $s$ and the triangle inequality give
\[
 \frac d m\le|\mu_t-\mu_s|
 \le|\mu_t-j|+|j-\mu_s|\le2|\mu_t-j|.
\]
Thus $|j-\mu_t|\ge d/(2m)$.
Hoeffding's inequality~\cite{Hoeffding63} for a sum $S$ of $d$
independent random variables with values in $[0,1]$ gives
\begin{equation}\label{eq:hoeffding}
 \Prob(|S-\E S|\ge u)\le2e^{-2u^2/d}.
\end{equation}
For $t\ne s$, by~\eqref{eq:bernoulli} the component $p_{L_t}$ is the law
of $1+S$, where $S$ is such a sum and $1+\E S=\mu_t$. Hence
\[
 p_{L_t}(j)\le\Prob\Bigl(|S-\E S|\ge\frac d{2m}\Bigr)
 \le2\exp\Bigl(-\frac2d\Bigl(\frac d{2m}\Bigr)^2\Bigr)
 =2\exp\Bigl(-\frac{n-1}{2m^2}\Bigr).
\]
Summing over the $m-1$ indices $t\ne s$ gives the upper bound
in~\eqref{eq:tails}; the lower bound is trivial.
\end{proof}

\subsection{A coefficient comparison for modal ties}
The remaining estimate is an exact comparison of adjacent coefficients,
used in Section~\ref{sec:strategy} to break ties between candidate labels. Coefficients outside a polynomial's support
are understood to be zero.

\begin{lemma}[Tie-breaking comparison]\label{lem:coefficients}
Let $D\ge1$ be an integer and $0<p<1$ arbitrary (in Section~\ref{sec:strategy}
the lemma is applied with $p=1/m$), and let
\[
 c_r=\coef{r}\prod_{\ell=1}^{D}(1-q_\ell+q_\ell x),
 \qquad q_\ell\in[p,1].
\]
For $1\le r\le D$,
\begin{equation}\label{eq:coefficientratio}
 (1-p)r c_r\ge p(D-r+1)c_{r-1}.
\end{equation}
In particular, $c_r\ge c_{r-1}$ if $r\le(D+1)p$.
If $q_\ell<1$ for all $\ell$ and $r<(D+1)p$, then $c_r>c_{r-1}$.
\end{lemma}

\begin{proof}
First suppose $q_\ell<1$ for all $\ell$, and put
$w_\ell=q_\ell/(1-q_\ell)\ge p/(1-p)$. Let
\[
 e_r=\sum_{\substack{A\subseteq\{1,\ldots,D\}\\|A|=r}}
     \prod_{\ell\in A}w_\ell
\]
be the elementary symmetric polynomial of degree $r$ in
$w_1,\ldots,w_D$, with $e_0=1$. Since
$\prod_\ell(1-q_\ell+q_\ell x)=\prod_\ell(1-q_\ell)\prod_\ell(1+w_\ell x)$,
we have $c_r=\prod_\ell(1-q_\ell)\,e_r$. Each $r$-subset arises exactly
$r$ times as $A\cup\{\ell\}$ with $|A|=r-1$ and $\ell\notin A$, and each
inner sum below has $D-r+1$ terms, each at least $p/(1-p)$; hence
\begin{align*}
 r e_r
 &=\sum_{\substack{A\subseteq\{1,\ldots,D\}\\|A|=r-1}}
       \left(\prod_{\ell\in A}w_\ell\right)
       \sum_{\ell\notin A}w_\ell\\
 &\ge(D-r+1)\frac{p}{1-p}e_{r-1}.
\end{align*}
Multiplication by $\prod_\ell(1-q_\ell)$ proves~\eqref{eq:coefficientratio}.
Parameters equal to $1$ follow by continuity. If every $q_\ell<1$, then
every $q_\ell\in(0,1)$ because $q_\ell\ge p>0$, so $c_{r-1}>0$;
the strict conclusion follows because
$p(D-r+1)>(1-p)r$ precisely when $r<(D+1)p$.
\end{proof}
\begin{remark}
$\lfloor(D+1)p\rfloor$ is a mode of $\Bin(D,p)$. The lemma asserts that
a Poisson-binomial distribution all of whose parameters are at least $p$
is still non-decreasing up to that point.
\end{remark}

\section{Proof of the optimal strategy}\label{sec:strategy}

We return to $m=2^k$, $p=1/m$, and the row-specific notation $N=i-1$,
$M=n-i$, $d=n-1$ of~\eqref{eq:auxiliary}. The integer $D$ of
Lemma~\ref{lem:coefficients} is specified separately in each application.

\begin{proof}[Proof of Theorem~\ref{thm:strategy}]
We first consider $i\le\lfloor(n+1)/2\rfloor$, so that $N\le M$.
Put
\[
 h=\lfloor(N+1)p\rfloor=\lfloor i/m\rfloor,
 \qquad \gamma_n=\frac1{m^2(n+1)^2}.
\]
The proposed guess is $h+1$, a mode of $p_{L_0}$. 
Choose $n(k)$ large enough that, for all $n\ge n(k)$,
\begin{equation}\label{eq:threshold}
 \varepsilon_n<\gamma_n,
 \qquad
 \max_j p_{L_t}(j)\le H_N-\gamma_n\quad(1\le t\le m-2)
\end{equation}
uniformly in the top-half position $i$. Here $\gamma_n$ is a common lower
bound for all the endpoint gaps defined after~\eqref{eq:endpoints}: since
$N+1\le n$ and $M+1\le n$, the right-hand sides of~\eqref{eq:pointgap}
(for $N$ and for $M$) and of~\eqref{eq:heightgap} are all at least
$\gamma_n$. Such a choice of $n(k)$ exists. The second condition follows
from~\eqref{eq:interiorgap}, because $\kappa_mn^{-1/2}\ge\gamma_n$ for large
$n$; the first holds eventually because $\varepsilon_n$ decays
exponentially while $\gamma_n$ is of exact order $n^{-2}$.
For $m=2$, only the first condition is needed.

At the proposed guess,
\[
 m a^{(k)}_{i,h+1}=\sum_{t=0}^{m-1}p_{L_t}(h+1)\ge H_N.
\]
\emph{Step 1: reduction to the endpoint modes.}
At any label $j$, choose its closest-mean component $s$ as in
Lemma~\ref{lem:tails}.
If $s$ is interior,~\eqref{eq:threshold} gives
$p_{L_s}(j)\le H_N-\gamma_n$.
If $s$ is an endpoint and $j$ is not a mode of that endpoint,
Lemma~\ref{lem:binomial}(i), together with $H_M\le H_N$, gives the same
bound. Finally, if $s=m-1$ and $H_M<H_N$, use
Lemma~\ref{lem:binomial}(ii). In all these cases,
\[
 m a^{(k)}_{i,j}\le p_{L_s}(j)+\varepsilon_n
 \le H_N-\gamma_n+\varepsilon_n<H_N.
\]
The first inequality is the upper bound in~\eqref{eq:tails}
of Lemma~\ref{lem:tails}.
Every maximizing label must therefore be a mode of the first component,
or a mode of the last component when $H_M=H_N$.

\emph{Step 2: the tie cases.} By Step~1, it remains to compare $h+1$ with
the other mode of the first component, if there is one, and, when
$H_M=H_N$, with the modes of the last component. There are three cases.

\emph{Case 2.1: two modes of the first component.}
The first component has two modes only when $ip$ is an integer, in which
case $h=ip$; their positions are then $h$ and $h+1$. For $t\ge1$, all
$d=n-1$ Bernoulli parameters in $F_{L_t}(x)/x$ are at least $p$.
Here $h=ip\le np=(d+1)p$. Moreover $1\le h\le d$: indeed $i/m\in\Z$ and
$i\ge1$ force $i\ge m$, so $h\ge1$, while
$h=i/m\le(n+1)/(2m)\le n-1=d$ for $n\ge2$.
Apply Lemma~\ref{lem:coefficients} to the polynomial $F_{L_t}(x)/x$, with
$D=d$ and $r=h$; its coefficients are
$c_r=\coef{r}\bigl(F_{L_t}(x)/x\bigr)=p_{L_t}(r+1)$, so the lemma gives
\[
 p_{L_t}(h+1)\ge p_{L_t}(h),\qquad 1\le t<m.
\]
The first component has equal probabilities at these labels. Summing
shows that $h+1$ is at least as good as $h$.
Figure~\ref{fig:mixture}
illustrates this case for $n=60$, $k=2$, $i=20$.

\emph{Case 2.2: equal endpoint heights, $n$ odd.}
If $M=N$, then $n$ is odd and $i=(n+1)/2$ (Figure~\ref{fig:peaks}(b)).
Equation~\eqref{eq:symmetry} makes the row invariant under
$j\mapsto n+1-j$. The endpoint modal sets are reflections of one another,
so the preceding comparison proves optimality of $h+1$ and its reflection.

\emph{Case 2.3: equal endpoint heights, $n$ even.}
Suppose $M>N$ and $H_M=H_N$. Lemma~\ref{lem:binomial}(ii) forces
$M=N+1=i$ and $ip=h\in\Z$. Hence $n=2i$.
The first component has modes $h,h+1$, and the last component has the
unique mode $n-h$. Define
\[
 A_t=1-tp+tp x,\qquad B_t=1-(t+1)p+(t+1)p x.
\]
Reflection of the paired component in~\eqref{eq:component} yields
\begin{align}
 F_{L_t}(x)-x^{n+1}F_{L_{m-1-t}}(1/x)
 &=xA_t^iB_t^{i-1}-xB_t^iA_t^{i-1}\notag\\
 &=p x(1-x)A_t^{i-1}B_t^{i-1}.\label{eq:reflectiondifference}
\end{align}
Let $c_{t,r}=\coef{r}A_t^{i-1}B_t^{i-1}$.
Taking the coefficient of $x^{h+1}$ and summing over $t$ gives
\begin{equation}\label{eq:tiebreak}
 a^{(k)}_{i,h+1}-a^{(k)}_{i,n-h}
 =\frac p m\sum_{t=0}^{m-1}(c_{t,h}-c_{t,h-1}).
\end{equation}
For $t=0$, the summand is zero: $A_0=1$, and $B_0^{i-1}$ is the
generating function of $\Bin(i-1,p)$, with modes $h-1,h$.
For $t\ge1$, apply Lemma~\ref{lem:coefficients} to the polynomial
$A_t^{i-1}B_t^{i-1}$, with $D=2i-2=d-1$ and $r=h$; its coefficients are
exactly $c_{t,r}$, and all its $D$ Bernoulli parameters are at least $p$.
Since $h=ip\le(2i-1)p=(D+1)p$, every remaining summand is nonnegative. Thus $h+1$ also beats or ties the
last-component mode.

These cases exhaust the possible maximizing labels in the top half.
Reflecting the resulting guesses by~\eqref{eq:symmetry} gives exactly
the second line of~\eqref{eq:strategy}, completing the proof.
\end{proof}

\begin{corollary}[Uniqueness for $k\ge2$]\label{cor:uniqueness}
Fix $k\ge2$. For all sufficiently large $n$, the maximizing label is
unique at every position except the middle position of an odd deck.
At that middle position the maximizing labels are exactly
\[
 \left\lfloor\frac{n+1}{2^{k+1}}\right\rfloor+1
 \quad\text{and}\quad
 n-\left\lfloor\frac{n+1}{2^{k+1}}\right\rfloor.
\]
\end{corollary}

\begin{proof}
The localization in Step~1 of the proof of Theorem~\ref{thm:strategy}
excludes all other labels strictly. In Case~2.1, use the component $t=1$:
its Bernoulli parameters are $p$ and $2p$, both in $(0,1)$ because
$m\ge4$. Since $h=ip<np$, the strict part of
Lemma~\ref{lem:coefficients} gives $p_{L_1}(h+1)>p_{L_1}(h)$. In Case~2.3,
the summand $t=1$ in~\eqref{eq:tiebreak} is likewise strictly positive,
since $h=ip<(2i-1)p$ and $i\ge m\ge4$. Thus all ties are resolved
strictly except for the exact reflection symmetry at the middle position
of an odd deck (Case~2.2), where the two reflected labels are distinct for
sufficiently large $n$.
\end{proof}

\begin{remark}[Ties]\label{rem:correction}
For $k\ge2$ and sufficiently large $n$, Corollary~\ref{cor:uniqueness}
shows that~\eqref{eq:strategy} is the only optimal strategy, apart from
the middle position of an odd deck. For $k=1$, the optimal
strategy is not unique, since there are no interior components to break
the ties.
\end{remark}

\section{The optimal expected score}\label{sec:expectation}

The score asymptotic requires only the height of the largest marginal
probability, rather than its exact location. The following estimate
therefore provides a proof independent of the tie analysis in
Section~\ref{sec:strategy}.

\begin{proposition}\label{prop:rowmaximum}
For fixed $k\ge1$ and all sufficiently large $n$, uniformly over
$1\le i\le\lfloor(n+1)/2\rfloor$,
\begin{equation}\label{eq:sandwich}
 \frac{H_{i-1}}m\le\max_j a^{(k)}_{i,j}
 \le\frac{H_{i-1}+\varepsilon_n}{m},
\end{equation}
where $\varepsilon_n$, defined in~\eqref{eq:epsilon}, bounds the total
contribution of all components other than the nearest one.
\end{proposition}

\begin{proof}
Evaluate~\eqref{eq:mixture} at any mode of the first component to obtain
the lower bound. For the upper bound, with $N=i-1$ and $M=n-i$ as
in~\eqref{eq:auxiliary}, Lemma~\ref{lem:binomial}(ii) gives $H_M\le H_N$, and Lemma~\ref{lem:interior} bounds all interior peaks
by $H_N$ once $n$ is sufficiently large. For $m=2$, there are only the
two endpoint components. Thus $p_{L_t}(j)\le H_N$ for all $t,j$.
Choose the closest-mean component at each $j$ and apply
Lemma~\ref{lem:tails}. Dividing by $m$ proves~\eqref{eq:sandwich}.
\end{proof}

\begin{proof}[Proof of Theorem~\ref{thm:expectation}]
Let $\ell=\lfloor n/2\rfloor$. Equations~\eqref{eq:objective}
and~\eqref{eq:symmetry}, followed by Proposition~\ref{prop:rowmaximum},
give
\begin{equation}\label{eq:sumheight}
 \E X^{G^*}=\frac2m\sum_{N=0}^{\ell-1}H_N+O_m(1).
\end{equation}
Indeed, the accumulated error from~\eqref{eq:sandwich} is
$O_m(n\varepsilon_n)=o(1)$, and when $n$ is odd the one unpaired
middle position contributes at most $1$.
Since $H_0=1$, Lemma~\ref{lem:binomial}(iii) implies
\begin{align*}
 \E X^{G^*}
 &=\frac{2}{m\sqrt{2\pi v}}
       \sum_{N=1}^{\ell-1}N^{-1/2}+O_m(1)\\
 &=\frac{4\sqrt\ell}{m\sqrt{2\pi v}}+O_m(1)\\
 &=\frac{2\sqrt n}{\sqrt{\pi(m-1)}}+O_m(1).
\end{align*}
Here the errors $O_m(N^{-3/2})$ are summable,
$\sum_{N=1}^{\ell-1}N^{-1/2}=2\sqrt\ell+O(1)$, and
$v=(m-1)/m^2$. Substituting $m=2^k$ proves~\eqref{eq:expectation}.
\end{proof}

\section{Further remarks}\label{sec:discussion}

The marginal representation from~\cite{KT25} separates two aspects of
optimal guessing. Peak heights determine the leading expected score,
whereas exact coefficient comparisons determine the strategy at ties.
Proposition~\ref{prop:rowmaximum} quantifies the first assertion with
an exponentially small error that is uniform in the position. Lemma~\ref{lem:coefficients} and the
reflection identity~\eqref{eq:reflectiondifference} address the second.

The threshold $n(k)$ in Theorem~\ref{thm:strategy} is not evaluated here.
The proof shows how it can be bounded: one needs a quantitative version
of the interior estimate~\eqref{eq:fourier} and a value of $n$ for which
$\varepsilon_n<\gamma_n$, with $\gamma_n$ as in Section~\ref{sec:strategy}. 
Numerically, the threshold appears to be small. Exact
computations show that the strategy~\eqref{eq:strategy} is optimal at
every position for $1\le n\le140$, $8\le n\le140$, $37\le n\le140$ and
$127\le n\le140$ when $k=1,2,3$ and $4$, respectively.

All results keep $k$ fixed. When $k$ grows with $n$, the number of
components increases and the separation estimate changes accordingly;
the arguments here do not provide a uniform transition to the
logarithmic-shuffle regime studied by Ciucu~\cite{Ciucu98}.
Another question is the distribution of $X^{G^*}$ for fixed $k\ge2$.
The single-position marginals suffice for its expectation, but higher
moments require joint probabilities. Extending the one-shuffle moment
and limit-law results of~\cite{KSTY23,KP24} therefore calls for additional
information about the dependence between positions.

The Maple programs accompanying this paper are available at
\url{https://thotsaporn.com/OptGuess.txt}.

\section*{Use of Generative-AI tools declaration}
All results and proof ideas in this paper are the authors' own.
Generative AI tools (OpenAI Codex and Anthropic Claude) were used to check
derivations, to verify computations independently, to suggest improvements
to the presentation of some arguments, and to edit the language of the
manuscript. The authors reviewed all content and take full responsibility
for the paper.

\section*{Conflict of interest}
The authors declare there is no conflict of interest.

\end{document}